\documentclass[11pt,reqno]{amsart}
\usepackage{latexsym}
\RequirePackage{amsmath}   \RequirePackage{amssymb}
\RequirePackage{amsthm}   \RequirePackage{epsfig}
\theoremstyle{plain}
\newtheorem{theorem}{Theorem}

\newtheorem*{corollary*}{Corollary}
\newtheorem{theoremO}{Theorem}
\newtheorem{lemmO}[theoremO]{Lemma}

\newtheorem*{conjecture*}{Conjecture}

\theoremstyle{definition}

\theoremstyle{remark}

\newcommand{\SC}{{\mathbb C}}  \newcommand{\SD}{{\mathbb D}}  
\newcommand {\SR}{{\mathbb R}}    
\newcommand{\al}{\alpha}    
  \newcommand{\ve}{\varepsilon}  \newcommand{\ze}{\zeta}
    \newcommand{\la}{\lambda}
\newcommand{\si}{\sigma}  \newcommand{\vp}{\varphi}  \newcommand{\om}{\omega}
  \newcommand{\Om}{\Omega}  
   
 \newcommand{\cF}{{\mathcal F}}

\newcommand{\be}{\begin{equation}}
\newcommand{\ee}{\end{equation}}
\newcommand{\bea}{\begin{eqnarray}}
\newcommand{\eea}{\end{eqnarray}}

\def\le{\left}
\def\ri{\right}

\begin{document}
\title{Ahlfors-Weill extensions for harmonic mappings}  

\author[I.~Efraimidis]{Iason Efraimidis}
\address{Department of Mathematics and Statistics, Texas Tech University, Box 41042, Lubbock, TX 79409, United States} \email{iason.efraimidis@ttu.edu}  

\author[R.~Hern\'andez]{Rodrigo Hern\'andez}
\address{Facultad de Ingenier\'{\i}a y Ciencias, Universidad Adolfo Ib\'a\~nez, Av. Padre Hurtado 750, Vi\~na del Mar, Chile} \email{rodrigo.hernandez@uai.cl}

\author[M.J.~Mart\'in]{Mar\'ia J.~Mart\'in} 
\address{Departamento de An\'alisis Matem\'atico, Universidad de La Laguna, 38200 San Crist\'obal de La Laguna, Santa Cruz de Tenerife, Spain} \email{maria.martin@ull.es} 


\maketitle

\begin{abstract}
We provide two new formulas for quasiconformal extension to $\overline{\SC}$ for harmonic mappings defined in the unit disk and having sufficiently small Schwarzian derivative. Both are generalizations of the Ahlfors-Weill extension for holomorphic functions. 
\end{abstract}

\section{Introduction}  

\subsection{Schwarzian derivative} Let $\SD$ be the unit disk and $f$ a locally univalent holomorphic function in $\SD$. The pre-Schwarzian and Schwarzian derivatives of $f$ are defined by
$$ 
Pf = (\log f')' = \frac{f''}{f'} \qquad \text{and} \qquad Sf =  (Pf)' -\tfrac{1}{2} (Pf)^2,  
$$
respectively, while the size of the Schwarzian derivative is measured by its norm 
\be \label{def-Schw-norm}
\|Sf\| \, = \, \sup_{z\in\SD} \, (1-|z|^2)^2 |Sf(z)|.
\ee 
These operators are closely related to criteria for univalence, as well as criteria for homeomorphic or quasiconformal extension. For example, Nehari's classical criterion states that if $\|Sf\|\leq 2$ then $f$ is injective in $\SD$. This was strengthened by Gehring and Pommerenke who proved that if $\|Sf\|\leq 2$ then the image $f(\SD)$ is a Jordan domain and $f$ admits a homeomorphic extension to $\overline{\SC}$, unless $f$ is M\"obius conjugate to the logarithm function $\ell(z)=\log\frac{1+z}{1-z}$, that is, $f=T\circ\ell\circ\tau$ for M\"obius transformations $T$ and $\tau$ with $\tau(\SD)=\SD$. Ahlfors and Weill \cite{AW62} showed that if $\|Sf\|\leq 2t$ for some $t<1$ then $f$ has a $\frac{1+t}{1-t}$-quasiconformal extension to $\overline{\SC}$. We refer the reader to the book of Lehto \cite{Leh} and the survey of Osgood \cite{Os} for more information on the Schwarzian derivative.

\subsection{Ahlfors-Weill extension} A remarkable feature of the Ahlfors-Weill theorem is that the extension is explicit. It is given by 
$$
F(z)  \,=\, \left\{
\begin{array}{rl}
f(z),   & \quad \text{if} \quad |z|\leq1,\\
E_f(1/\overline{z}), & \quad \text{if} \quad |z|>1, 
\end{array} \right.
$$
where
\be \label{ext-E-f}
E_f(\ze)  = f(\ze)+  \frac{(1-|\ze|^2)f'(\ze)}{\overline{\ze}-\tfrac{1}{2}(1-|\ze|^2) Pf(\ze)},  \qquad \ze \in\SD. 
\ee
The Beltrami coefficient $\mu_F=F_{\overline{z}} / F_z$ of the extension can be computed as 
$$
\mu_F(z) = - \frac{1}{2} \left( \frac{\ze}{\overline{\ze}} \right)^2 (1-|\ze|^2)^2 Sf(\ze), \qquad |z|>1, \; \ze=1/\overline{z},  
$$
so that the condition $\|Sf\|\leq 2t$ implies that $|\mu_F|\leq t <1$. Hence its Jacobian $J_F=|F_z|^2-|F_{\overline{z}}|^2$ is positive and, therefore, $F$ is locally homeomorhic in $|z|>1$. After showing that $f$ is continuous up to the boundary $\partial\SD$ and it matches the extension there, Ahlfors and Weill deduce that $F$ is locally homeomorhic in the Riemann sphere $\overline{\SC}$. From there, the following well-known topological proposition concludes their proof.

\begin{lemmO}\label{lem-top-monodromy}
A locally homeomorhic mapping $F:\overline{\SC}\longrightarrow \overline{\SC}$ is a homeomorphism. 
\end{lemmO}

See \cite[p.~23]{Ah74} or \cite[\S15C]{AlSa} for this. The standard proof involves first showing that $F(\overline{\SC})= \overline{\SC}$, then that $F$ is a covering map and finally applying the Monodromy theorem. See \cite{Ka94} for a different approach.

Most extension results involving the Schwarzian derivative in the literature rely on  the use of this lemma; see, for example, \cite{Ah74, AW62, Ch20, CDO18}. An exception is a result of Chuaqui and Osgood \cite{CO94} who showed that under the weaker hypothesis $\|Sf\|\leq 2$ the Ahlfors-Weill extension is, still, homeomorphic if $f$ is not M\"obius conjugate to the logarithm function. To prove this they studied the critical points of the density $\la_\Om$ of the hyperbolic metric in $\Om=f(\SD)$, defined by
$$
\la_\Om\big(f(z)\big) |f'(z)| = \la_\SD(z) = \frac{1}{1-|z|^2}, \qquad z\in\SD, 
$$
and used the fact that 
\be \label{form-ext-hyper}
E_f(\ze)  = w + \frac{1}{\partial_w \log \la_\Om(w)}, \qquad w=f(\ze), \; \ze \in \SD.  
\ee 
A geometric interpretation of the extension follows from this formula. Simply note that $\partial_z u = \frac{1}{2}\overline{\nabla u}$ for any function $u$, so that the complex numbers $\big(\partial_w \log \la_\Om(w)\big)^{-1}$ and $ \nabla \log \la_\Om(w)$ have the same argument. Hence, the extension moves away from $w$ in the direction of maximal growth of $\la_\Om$. Minda \cite{Mi97} gave further geometric insight for formula \eqref{form-ext-hyper}.

\subsection{Harmonic mappings} A complex-valued harmonic mapping $f$ in $\SD$ has a canonical decomposition $f=h+\overline{g}$, where $h$ and $g$ are analytic in $\SD$. By Lewy's theorem, $f$ is locally univalent if and only if its Jacobian $J_f=|h'|^2-|g'|^2$ does not vanish. Also, $f$ is orientation-preserving if its dilatation $\omega=g'/h'$ satisfies $|\omega|<1$ in $\SD$.  See Duren's book \cite{Du} for more on the theory of planar harmonic mappings.

A generalization of the Schwarzian derivative to harmonic mappings was first given by Chuaqui, Duren and Osgood \cite{CDO03} and later by the second and third authors \cite{HM15} of the present paper. We will follow the latter, which does not involve the Weierstarss-Enneper lift to a minimal surface and so it fits better to considerations within the complex plane. Hence, the Schwarzian derivative of a locally univalent harmonic mapping $f$ is defined by 
\be \label{def-Sch-HM}
S_f \, = \, (P_f)_z - \tfrac{1}{2} (P_f)^2, 
\ee
where $ P_f =(\log J_f)_z$ is the pre-Schwarzian derivative. These are equivalent to 
\be \label{pre-Sch-HM}
P_f \, = \, Ph - \frac{\overline{\om}\om'}{1-|\om|^2}
\ee
and
\be \label{Sch-HM}
S_f \, = \, Sh + \frac{\overline{\om}}{1-|\om|^2} \left( \frac{h''}{h'} \om'-\om'' \right) - \frac{3}{2} \left( \frac{\om' \overline{\om}}{1-|\om|^2} \right)^2. 
\ee
Note that we are using the notation $Pf$ and $Sf$ when we know that the mapping $f$ is holomorphic and the notation $P_f$ and $S_f$, with the mapping $f$ as a subscript, in the more general setting of harmonic mappings. The Schwarzian norm of a harmonic mapping $f$ is defined exactly as in \eqref{def-Schw-norm}. See \cite{HM15} for more about the basic properties of these operators.

The second and third authors showed in \cite{HM15} that $S_f\equiv0$ implies that $f$ is a harmonic M\"obius transformation, \emph{i.e.}, post-composition of a M\"obius transformation with an affine map. Moreover, the same authors proved in \cite{HM15-2} that there exists a constant $c>0$ such that if $\|S_f\|\leq c$ then $f$ is univalent in $\SD$. It is an interesting open problem to find the larger such constant $c$ or, at least, give some positive lower bound for it. A criterion for the existence of a quasiconformal extension was also given in \cite{HM15-2}. These results have been extended by the first author in \cite{Ef1,Ef2} to more general domains.

Following their definition \cite{CDO03} of the Schwarzian derivative, Chuaqui, Duren and Osgood found in \cite{CDO18} a quasiconformal extension to the three-dimensional space $\SR^3$ for the Weierstarss-Enneper lift of a harmonic mapping in $\SD$. They showed that under an additional assumption the image of $\SC$  is locally a graph and, therefore, projecting it to $\SC$ results in a planar quasiregular mapping. In view of Lemma~\ref{lem-top-monodromy} this mapping is actually a quasiconformal extension of the initial harmonic mapping. This result reduces to the original Ahlfors-Weill theorem when the given mapping in $\SD$ is analytic and, moreover, it gives the correct constants: the bound $2t$ on the Schwarzian yields a $\frac{1+t}{1-t}$-quasiconformal extension to $\overline{\SC}$. With the same approach, Chuaqui \cite{Ch20} provided many more extensions. 

Note that harmonic mappings admitting a Weierstarss-Enneper lift to a minimal surface require that their dilatation is the square of a meromorphic function. The purpose of this paper is to remove this hypothesis and give two formulas for quasiconformal extension to $\overline{\SC}$ by using the Schwarzian derivative \eqref{Sch-HM} and techniques entirely within the complex plane.

\subsection{Main result} We may now state our main result.

\begin{theorem} \label{main-thm}
Let $d\in[0,1)$ and $f=h+\overline{g}$ be a locally univalent harmonic mapping in $\SD$ whose dilatation $\om$ satisfies $|\om(z)|\leq d$ for all $z\in\SD$. Let also $\tau$ be either 
\begin{itemize}
\item[(A)] the pre-Schwarzian $Ph$ or
\item[(B)] the pre-Schwarzian $P_f$.
\end{itemize}
Then, for either case {\rm(A)} or {\rm (B)}, there exists a constant $c=c(d)>0$ such that if $\|S_f\|\leq c$ then the mapping 
\be \label{ext-disk-out}
F(z)  \,=\, \left\{
\begin{array}{rl}
f(z),   & \quad \text{if} \quad |z|\leq1,\\
E_f(1/\overline{z}), & \quad \text{if} \quad |z|>1, 
\end{array} \right.
\ee
where
$$
E_f(\ze)  = f(\ze)+ \Phi(\ze) + \overline{\om(\ze)} \overline{\Phi(\ze)}, \qquad \ze \in\SD, 
$$
and 
\be \label{Phi}
\Phi(\ze) = \frac{(1-|\ze|^2)h'(\ze)}{\overline{\ze}-\tfrac{1}{2}(1-|\ze|^2) \tau(\ze)},  \qquad \ze \in\SD, 
\ee
is a quasiconformal extension of $f$ to $\overline{\SC}$.
\end{theorem}

We note that the homeomorphic extension of $f$ to the boundary $\partial \SD$ was given in \cite{Ef1} under the assumption that  $\|S_f\|$ is sufficiently small and under no assumption on the dilatation $\om$. 

The extension for case (A) is constructed via best M\"obius approximation in Subsection~\ref{sub-sect-best-Mob}, while the construction for case (B) is given in Subsection~\ref{sub-sect-conf-fact} and it involves choosing a different conformal factor than the one chosen in \cite{CDO18}. A connection between the two extensions is given in Subsection~\ref{sub-sect-prop-conn}. The proof of the theorem, given in Section~\ref{Sect-proof-A} for case (A) and Section~\ref{Sect-proof-B} for case (B), relies on approximating $f$ by its dilation $f(rz), r<1$, showing that the Beltrami coefficient of the extension is bounded away from 1 and invoking Lemma~\ref{lem-top-monodromy}. The desired quasiconformal extension is then harvested as the limit for $r\to1$.

\section{Construction and properties of the extensions}  

\subsection{Best M\"obius approximation and the extension (A)} \label{sub-sect-best-Mob} Martio and Sarvas~\cite{MS79} were the first to use best M\"obius approximation in order to prove criteria for injectivity. Moreover, the observation that the Ahlfors-Weill extension can also be obtained via best M\"obius approximation seems to be well known; see \cite{Mi97} or \cite[\S3.3]{Os}. We briefly describe it here and then extend this idea to harmonic mappings.

If $f$ is analytic at a point $\ze$ then there is a unique M\"obius transformation $M=M(f,\ze)$ that agrees with $f$ up to second order at $\ze$, that is, it satisfies $M^{(k)}(\ze)=f^{(k)}(\ze)$, for $k=0,1,2$. It is easy to see that if $g$ is analytic at the origin and normalized by $g(0)=g'(0)-1=0$ then its best M\"obius approximation at the origin is 
$$
M(g,0)(z) = \frac{z}{1-\tfrac{1}{2}g''(0)z}. 
$$
If now $f:\SD\to\SC$ is locally univalent and $\ze$ is a point in $\SD$ then 
$$
g(z) = \frac{f(z+\ze)-f(\ze)}{f'(\ze)} 
$$
is analytic at the origin and normalized as before. (Alternatively, we may normalize $f$ by pre-composing with an appropriate disk automorphism.) Hence, the best M\"obius approximation of $f$ at $\ze$ is given by
$$
M(f,\ze)(z) = f(\ze)+ f'(\ze) M(g,0)(z-\ze) = f(\ze)+ \frac{(z-\ze)f'(\ze)}{1-\tfrac{1}{2}(z-\ze)Pf(\ze)}. 
$$
Now the Ahlfors-Weill extension \eqref{ext-E-f} can be obtained by setting 
$$
E_f(\ze)= M(f,1/\overline{z})(z), \qquad \ze=1/\overline{z}, \;\;  |z|>1. 
$$
Note that for a fixed $|z|>1$, the points $\ze$ and $z$ are symmetric with respect to the unit circle $\partial \SD$ and, therefore, the M\"obius map $M(f,\ze)$ sends them to the points $f(\ze)$ and $F(z)$, respectively, which are symmetric with respect to the circle $M(f,\ze)(\partial\SD)$.

A harmonic M\"obius  transformation is a mapping $M=T+a \overline{T}$, where $a\in\SD$ is a constant and $T$ is a (holomorphic) M\"obius  transformation; see \cite{CDO03}. If $f$ is harmonic at a point $\ze$ then the harmonic M\"obius  transformation that best approximates $f$ at $\ze$ is the unique mapping $M=M(f,\ze)$ that agrees with $f$ at the value, the first two analytic derivatives and the first anti-analytic derivative, that is, it satisfies 
$$
M(\ze)=f(\ze), \qquad M_z(\ze)=f_z(\ze),  \qquad M_{zz}(\ze)=f_{zz}(\ze),  \qquad M_{\overline{z}}(\ze)=f_{\overline{z}}(\ze). 
$$
If $f=h+\overline{g}$ is normalized so that $h(0)=g(0)=0$ and $h'(0)=1$ then its best harmonic M\"obius approximation at the origin is given by 
\be \label{best-Mob-normal-harm}
M(f,0)(z) = \frac{z}{1-\tfrac{1}{2}h''(0)z} +\overline{g'(0)} \overline{\le( \frac{z}{1-\tfrac{1}{2}h''(0)z} \ri)}. 
 \ee
Moreover, if $f=h+\overline{g}$ is an arbitrary locally univalent harmonic mapping in $\SD$ and $\ze$ is a point in $\SD$ then we may normalize it as before by setting 
\be \label{normalized-harmonic}
F(z) = \frac{f(z+\ze)-f(\ze)}{h'(\ze)}. 
\ee
Writing $F=H+\overline{G}$ we see that $H''(0)=Ph(\ze)$ and $G'(0) = g'(\ze) / \overline{h'(\ze)}$. Hence, the best harmonic M\"obius approximation of $f$ at $\ze$ is given by 
\begin{align*}
M(f,\ze)(z) & = f(\ze)+ h'(\ze) M(F,0)(z-\ze)  \\
& = f(\ze)+  \frac{(z-\ze)h'(\ze)}{1-\tfrac{1}{2}(z-\ze)Ph(\ze)} + \overline{\om(\ze)} \overline{\le( \frac{(z-\ze)h'(\ze)}{1-\tfrac{1}{2}(z-\ze)Ph(\ze)} \ri)}. 
\end{align*}
For $|z|>1$, setting  
$$
E_f(\ze)=M(f,1/\overline{z})(z)  = f(\ze)+ \Phi(\ze) + \overline{\om(\ze)} \overline{\Phi(\ze)}, \qquad \ze=1/\overline{z}, 
$$
where
$$
\Phi(\ze) = \frac{(1-|\ze|^2)h'(\ze)}{\overline{\ze}-\tfrac{1}{2}(1-|\ze|^2)Ph(\ze)} \qquad \ze \in\SD,
$$
we get the extension \eqref{ext-disk-out} of $f$ to $\overline{\SC}$ in the case (A) of Theorem~\ref{main-thm}. Now the symmetric with respect to $\partial \SD$ points $\ze$ and $z$ are mapped by $M(f,\ze)$ to the points $f(\ze)$ and $F(z)$, respectively, which are symmetric with respect to the ellipse $M(f,\ze)(\partial\SD)$. Here for two points to be symmetric with respect to an ellipse we understand that they lie on the same ray emanating from the center of the ellipse (center of inversion) and that the product of their distances to the center equals the square of the distance from the center to the intersection of the ray on which they lie with the ellipse (see \cite{Ch65}).

\subsection{A different conformal factor and the extension (B)} \label{sub-sect-conf-fact} Chuaqui, Duren and Osgood \cite{CDO18} considered harmonic mappings $f=h+\overline{g}$ in $\SD$ whose dilatation is the square of an analytic function and, under certain assumptions, extended them with \eqref{ext-disk-out} and the formula 
$$
E_f(\ze)=  f(\ze)+ \frac{(1-|\ze|^2)h'(\ze)}{\overline{\ze}- (1-|\ze|^2) \partial_z \si(\ze)} + \frac{(1-|\ze|^2) \overline{g'(\ze)}}{ \ze- (1-|\ze|^2) \partial_{\overline{z}} \si(\ze)}, \quad \ze \in\SD. 
$$
Here $e^\si = |h'|+|g'|$ is the conformal factor of the metric induced by the Weierstarss-Enneper lift on the minimal surface. If instead of this conformal factor we use the Jacobian and take $\si=\log\sqrt{J_f}$ then we find that 
$$
\partial_z \si = \frac{1}{2} \le(Ph - \frac{\om'\overline{\om}}{1-|\om|^2}\ri) = \frac{1}{2} P_f \qquad \text{and} \qquad \partial_{\overline{z}}\si =  \frac{1}{2} \overline{P_f}.  
$$
Hence, for this choice of $\si$ the above becomes $E_f = f+ \Phi + \overline{\om \Phi}$, where
$$
\Phi(\ze) = \frac{(1-|\ze|^2)h'(\ze)}{\overline{\ze}-\tfrac{1}{2}(1-|\ze|^2)P_f(\ze)}, \qquad \ze\in\SD, 
$$
and therefore gives us the extension \eqref{ext-disk-out} of $f$ to $\overline{\SC}$ in the case (B) of Theorem~\ref{main-thm}.

Another way of arriving at this formula can be seen as follows. We consider an arbitrary locally injective mapping $f=h+\overline{g}$ in $\SD$ and a point $\ze\in\SD$ and, after normalizing as in \eqref{normalized-harmonic} to get $F$, we normalize further by post-composing with an affine map in order to get a mapping whose first anti-analytic coefficient is zero. Hence, we write 
$$
\widehat{f} = \frac{ F-\overline{G'(0)F} }{ 1-|G'(0)|^2 } = \widehat{h} + \overline{ \widehat{g}  }
$$
and easily see that $\widehat{g}'(0)=0$. Now, in view of \eqref{best-Mob-normal-harm}, we have that the best harmonic M\"obius approximation of $\widehat{f}$ at the origin is 
$$
M(\widehat{f},0)(z) =\frac{z}{1-\tfrac{1}{2}\widehat{h}''(0)z}, 
$$
coinciding with the approximation for $\widehat{h}$, for which we may also compute that $\widehat{h}''(0) = P_f(\ze)$. Inverting the above affine transformation we get $F= \widehat{f} + \overline{ G'(0)\widehat{f}}$ which motivates us to write the harmonic M\"obius transformation 
\be \label{non-best-Mob-approx}
M(\widehat{f},0) + \overline{ G'(0) M(\widehat{f},0)}.  
\ee
Now, this is \emph{not} the best harmonic M\"obius approximation $M(F,0)$ since the second analytic derivative $F_{zz}$ involves the second anti-analytic derivative $\widehat{f}_{\overline{zz}}$ and this information is encoded in $M(F,0)$, but not in \eqref{non-best-Mob-approx}. In short, the concept of best harmonic M\"obius approximation is not affine invariant. However, we may use \eqref{non-best-Mob-approx} as a (non-best) approximation of $F$ at the origin. From there, dismantling the normalization \eqref{normalized-harmonic} leads us to a (non-best) harmonic M\"obius approximation of $f$ at $\ze$ given by 
$$
R(f,\ze)(z) = f(\ze)+  \frac{(z-\ze)h'(\ze)}{1-\tfrac{1}{2}(z-\ze)P_f(\ze)} + \overline{\om(\ze)} \overline{\le( \frac{(z-\ze)h'(\ze)}{1-\tfrac{1}{2}(z-\ze)P_f(\ze)} \ri)}. 
$$
The extension in the case (B) of Theorem~\ref{main-thm} can now be obtained by setting $R(f,1/\overline{z})(z)$ for $|z|>1$. Even though unconventional, this argument permits the same geometric conclusions for case (B) as in case (A), \emph{i.e.}, that points $\ze$ and $z$ which are symmetric with respect to the unit circle $\partial \SD$ are mapped by $R(f,\ze)$ to the points $f(\ze)$ and $F(z)$, respectively, which are symmetric with respect to the ellipse $R(f,\ze)(\partial\SD)$.

\subsection{Properties}  \label{sub-sect-prop-conn} Straightforward calculations can show that both extensions (A) and (B) satisfy 
$$
E_{af+b} = aE_f+b, \qquad a,b\in\SC, 
$$
as well as 
$$
E_{f\circ\vp_\al} = E_f\circ \vp_\al, \qquad \al\in\SD, 
$$ 
where 
\be \label{disk-auto}
\vp_\al(z) = \frac{z-\al}{1-\overline{\al}z}, \qquad z \in\SD, 
\ee 
is a disk automorphism. Moreover, of the two extensions only (B) satisfies the affine invariance 
$$
E_{f+a\overline{f}}=E_f + a\overline{E_f}, \qquad a\in\SD. 
$$

To see a connection between the extensions (A) and (B) let us momentarily denote by $\Phi_f^j$ and $E_f^j$, for $j={\rm A, B}$, the ingredients of the extensions (A) and (B), respectively. We fix $|z|>1$, take $\ze=1/\overline{z}$ and consider the affine transformation $\widehat{f}=f-\overline{\om(\ze)f}$. We write $\widehat{f}=\widehat{h}+\overline{\widehat{g}}$, where
$$
\widehat{h} = h -\overline{\om(\ze)}g \qquad \text{and} \qquad \widehat{g}=g-\om(\ze)h, 
$$
and compute the dilatation of $\widehat{f}$ as 
$$
\widehat{\om} = \frac{ \om-\overline{\om(\ze)} }{ 1-\overline{\om(\ze)}\om }. 
$$
Evidently $\widehat{\om}(\ze) =0$. Also, 
$$
P\widehat{h} = Ph - \frac{ \overline{\om(\ze)} \om' }{ 1-\overline{\om(\ze)}\om },
$$
so, in view of \eqref{pre-Sch-HM}, we have that $P\widehat{h}(\ze)=P_f(\ze)$. Now, it is easy to see that 
$$
\Phi^{\rm A}_{\widehat{f}}(\ze) = (1-|\om(\ze)|^2)\Phi^{\rm B}_{f}(\ze), 
$$
which implies that
\begin{align*}
E^{\rm A}_{\widehat{f}}(\ze) & = f(\ze) - \overline{\om(\ze)f(\ze)} + (1-|\om(\ze)|^2)\Phi^{\rm B}_{f}(\ze) \\
& =  f(\ze) + \Phi^{\rm B}_{f}(\ze) + \overline{\om(\ze) \Phi^{\rm B}_{f}(\ze)} - \overline{\om(\ze)} \overline{ \left( f(\ze) + \Phi^{\rm B}_{f}(\ze) + \overline{\om(\ze) \Phi^{\rm B}_{f}(\ze)} \right) }\\
& = E^{\rm B}_f(\ze) -\overline{ \om(\ze) E^{\rm B}_f(\ze) }, 
\end{align*}
that is, an affine transformation (which depends on the point of evaluation) of the extension $E^{\rm B}_f$. This can also be inverted in order to give
$$
E^{\rm B}_f(\ze) = \frac{ E^{\rm A}_{\widehat{f}}(\ze) + \overline{ \om(\ze) E^{\rm A}_{\widehat{f}}(\ze) } }{ 1-|\om(\ze)|^2 }, \qquad \ze \in\SD. 
$$
However, this will not be used in our proofs.

\section{Preliminaries}  

\subsection{Bounded Schwarzian derivative} According to Pommerenke's \cite{Po64} well-known theorem, if $h$ is analytic and locally univalent in $\SD$ then 
\be  \label{Pom-lem}
(1-|z|^2)\left| \frac{h''(z)}{h'(z)}  -\frac{2 \, \bar{z}}{1-|z|^2} \right| \, \leq \, 2 \sqrt{ 1+\frac{\|Sh\| }{2}}, \qquad z\in\SD. 
\ee

We denote by $\cF_t$ the set of all sense-preserving harmonic mappings $f=h+\overline{g}$ in $\SD$ which satisfy $\|S_f\|\leq t$, for $t>0$, and are normalized by $h(0)=g(0)=0$ and $h'(0)=1$. This constitutes a normal family (see \cite[\S2.2]{Ef1}).

\subsection{The hyperbolic derivative} If $\om$ is an analytic self-mapping of $\SD$ then its hyperbolic derivative is given by 
\be \label{def-hyp-derivative}
\om^*(z) \, =\, \frac{(1-|z|^2)\om'(z)}{1-|\om(z)|^2}, \qquad z\in\SD, 
\ee
and the quantity $\|\om^*\| = \sup_{z\in\SD} |\om^*(z)|$ is called the hyperbolic norm of $\om$. In view of the Schwarz-Pick lemma we have that $\|\om^*\|\leq1$. Moreover, it has been shown in \cite{GZ} that 
\be \label{ineq-hyp-norm}
\frac{(1-|z|^2)^2 |\om''(z)|}{1-|\om(z)|^2} \, \leq \, C \|\om^*\|, \qquad z\in\SD, 
\ee
for some constant $C>0$.  

Let 
$$  
R_t \,=\, \sup \|\om^*\|, 
$$  
where the supremum is taken over all admissible dilatations for mappings in $\cF_t$. A compactness argument for a suitably normalized subclass of $\cF_t$ was used in \cite{HM15-2} to show that $R_t\to0$ as $t\to0$. The exact asymptotics of this convergence are not known, although in \cite[\S5.1]{CHM} it was shown that $R_t > \root 4 \of{2t/3}$, for $t\in(0,3/2)$.

In view of \eqref{Sch-HM}, for any mapping $f=h+\overline{g}$ we have that 
$$
|Sh| \, \leq \, |S_f| + \frac{|\om'|}{1-|\om|^2} \left| \frac{h''}{h'}\right|  + \frac{|\om''|}{1-|\om|^2} + \frac{3}{2} \left( \frac{|\om'|}{1-|\om|^2} \right)^2. 
$$
If $f  \in \cF_t$ then we may use \eqref{Pom-lem} and \eqref{ineq-hyp-norm} in order to obtain 
\be \label{Sh-small}
\| Sh \| \leq t + C R_t
\ee
for some $C>0$. This calculation first appeared in \cite{HM15-2}.

\section{Proof of Theorem~\ref{main-thm}(A): case $\tau=Ph$ }   \label{Sect-proof-A}
Let $f\in\cF_t$, for $t>0$, with $|\om|\leq d$. We consider the dilation $f_r(z)=f(rz)$, for $r<1$, with representation $f_r=h_r+\overline{g_r}$ and dilatation $\om_r(z)=\om(rz)$. It is clear from \eqref{Phi} that, since both $h_r'$ and $Ph_r$ are bounded, the mapping $E_{f_r}$ matches $f_r$ on the boundary $\partial \SD$. Also, the denominator of $\Phi$ is continuous in $\SD$. Therefore, the extension $F_r$ of $f_r$ is continuous in $\overline\SC$ with respect to the spherical metric. Both mappings $f_r$ and $E_{f_r}$ have derivatives (of all orders) in an open set containing $\overline{\SD}$, hence so does $F_r$ in $\overline\SC$. Since $|\om_r|\leq d <1$, it suffices to show that the Beltrami coefficient of $F_r$ in $|z|>1$ is bounded away from 1, in order to conclude that $F_r$ is locally homeomorphic in $\overline\SC$. If this is proved then Lemma~\ref{lem-top-monodromy} will show that $F_r$ is actually homeomorphic, hence quasiconformal, in $\overline\SC$. In view of Theorem 5.3 in \cite[Ch.II, \S 5.4]{LV}, letting $r\to1$  the limit mapping $F$ is either constant, or takes two values, or is quasiconformal. The first two cases are discarded by the normalization at the origin, so that the mapping $F$ in the statement of the theorem is quasiconformal in $\overline{\SC}$.

We will compute the Beltrami coefficient of the extension for the initial mapping $f$. The computation can clearly then be applied to its dilation $f_r$. We have that
$$
\Phi_\ze = \frac{h' \le[ \tfrac{1}{2} (1-|\ze|^2)^2 Sh - \tfrac{1}{4} (1-|\ze|^2)^2 \le(Ph\ri)^2 + (1-|\ze|^2) \, \overline{\ze} \, Ph - \overline{\ze}^2  \ri] }{\le[\overline{\ze}-\tfrac{1}{2}(1-|\ze|^2)Ph\ri]^2}
$$
and 
\be \label{Phi-zed-bar}
\Phi_{\overline{\ze}} = \frac{- h'}{\le[\overline{\ze}-\tfrac{1}{2}(1-|\ze|^2) Ph \ri]^2}. 
\ee
Now, for $|z|>1$ and $\ze=1/\overline{z}\in\SD$,  we have 
\be \label{F-zed}
F_z = -\frac{1}{z^2} \le( \overline{g'} +\Phi_{\overline{\ze} } + \overline{\om' \Phi }  + \overline{\om \Phi_\ze}  \ri)
\ee
and
\be \label{F-zed-bar}
F_{\overline{z}} = -\frac{1}{\overline{z}^2} \le( h' +\Phi_\ze + \overline{\om \Phi_{\overline{\ze}} }  \ri). 
\ee
Hence the Beltrami coefficient of the extension is
$$
\mu_F = \frac{F_{\overline{z}}}{F_z } = \le(\frac{z}{\overline{z}}\ri)^2 \frac{  h' +\Phi_\ze + \overline{\om \Phi_{\overline{\ze}} } }{   \overline{\om h'} +\Phi_{\overline{\ze} } + \overline{\om' \Phi }  + \overline{\om \Phi_\ze}    }. 
$$
We set 
$$
D = \overline{\ze}-\tfrac{1}{2}(1-|\ze|^2) Ph(\ze). 
$$
for the denominator of $\Phi$. Let also 
\be \label{lambda-T}
\la = \frac{\overline{h'}\, D^2}{h' \, \overline{D^2}} \, \in \partial\SD.
\ee 
Multiplying the numerator of $\mu$ by $D^2/h'$, the denominator by $\overline{D^2}/\overline{h'}$ and using \eqref{Phi-zed-bar} we obtain 
$$
\le(\frac{\overline{z}}{z}\ri)^2 \mu_F = \overline{\la} \, \frac{ D^2+\tfrac{D^2}{h'}\Phi_\ze -\la \overline{\om} }{  \overline{\om} \overline{ \le( D^2+\tfrac{D^2}{h'}\Phi_\ze \ri)} -\overline{\la} + \overline{\om' \tfrac{D^2}{h'}\Phi}  }. 
$$
It is straightforward to see that  
$$
D^2+\frac{D^2}{h'}\Phi_\ze = \frac{1}{2}(1-|\ze|^2)^2 \,Sh. 
$$
Hence, we get that 
\begin{align} \label{A-mu-fraction}
\le(\frac{\overline{z}}{z}\ri)^2 \mu_F & = \overline{\la} \, \frac{ \tfrac{1}{2}(1-|\ze|^2)^2 Sh -\la \overline{\om}   }{  \overline{\om} \tfrac{1}{2}(1-|\ze|^2)^2 \overline{Sh} -\overline{\la} + \overline{\om'}(1-|\ze|^2)\overline{D}  } \nonumber \\ 
& =  \frac{ \la \overline{\om} -\tfrac{1}{2}(1-|\ze|^2)^2 Sh  }{  1-  \la\overline{\om} \tfrac{1}{2}(1-|\ze|^2)^2 \overline{Sh}  - \la \overline{\om'}(1-|\ze|^2)\overline{D}  }.
\end{align}
We set 
$$
\al = \la \overline{\om} \qquad \text{and} \qquad \beta = \tfrac{1}{2}(1-|\ze|^2)^2 Sh
$$
and note that by the hypothesis and by \eqref{Sh-small} we have that 
$$
|\al| \leq d <1 \qquad \text{and} \qquad |\beta| \leq  \tfrac{1}{2}(t+C R_t)<1, 
$$
if $t$ is assumed to be sufficiently small. Dividing the numerator and the denominator of \eqref{A-mu-fraction} by $1- \overline{\beta}\al $ and using the expression \eqref{def-hyp-derivative} for the hyperbolic derivative we get that 
$$
\le( \frac{\overline{z}}{z} \ri)^2 \mu_F =\frac{ \vp_\beta(\al) }{ 1 - \frac{ \la \overline{\om^*}(1-|\om|^2)\overline{D} }{ 1- \overline{\beta}\al  } }, 
$$
where $\vp_\beta(\al)$ is a disk automorphism as in \eqref{disk-auto}. Therefore, we have that 
$$
|\mu_F| =  \frac{ |\vp_\beta(\al)| }{  \le|1   - \frac{ \la \overline{\om^*}(1-|\om|^2)\overline{D} }{ 1- \overline{\beta}\al  }\ri| } \leq  \frac{ |\vp_\beta(\al)| }{  1   - \frac{ |\om^*| |D| }{ |1- \overline{\beta}\al|  } }. 
$$
Now, it obviously holds that
$$
 |1- \overline{\beta}\al|  \geq 1- |\al\beta| \geq 1- \tfrac{d}{2}(t+C R_t) 
$$
and, by the Schwarz-Pick lemma, we have that 
$$
|\vp_\beta(\al)|  \leq \frac{|\al| + |\beta|}{1+ |\al\beta|} \leq \frac{d + \tfrac{1}{2}(t+C R_t)}{1+ \tfrac{d}{2}(t+C R_t)}. 
$$
Also, recalling that $|\om^*|\leq R_t$ and that \eqref{Pom-lem} and \eqref{Sh-small} yield 
$$
|D|\leq \sqrt{ 1+\frac{1}{2} \|Sh\| } \leq \sqrt{ 1+\frac{1}{2}(t+C R_t ) }, 
$$
we obtain  
$$
|\mu_F|  \leq  \frac{ \frac{d + \tfrac{1}{2}(t+C R_t)}{1+ \tfrac{d}{2}(t+C R_t)} }{ 1 - \frac{ R_t }{ 1- \tfrac{d}{2}(t+C R_t)   } \sqrt{ 1+\frac{1}{2} (t+C R_t) } }.  
$$
Since the bound for $|\mu_F|$ is a constant that tends to $d$ as $t\to0$, we have that for each $0<\ve<1-d$ there exists $t$ sufficiently small so that $\sup_{|z|>1}|\mu(z)|\leq d+ \ve<1$.

\section{Proof of Theorem~\ref{main-thm}(B): case $\tau=P_f$}  \label{Sect-proof-B}
The line of reasoning here follows the previous proof \emph{in toto}: we consider the dilation $f_r(z)=f(rz)$, for $r<1$, see that its extension $F_r$ is continuous in $\overline\SC$ with respect to the spherical metric, bound the Beltrami coefficient of $F_r$ in $|z|>1$, use Lemma~\ref{lem-top-monodromy} and pass to the limit for $r\to1$.

We set 
$$
D = \overline{\ze}-\tfrac{1}{2}(1-|\ze|^2) P_f(\ze)
$$
for the denominator of $\Phi$ and compute 
$$
\Phi_{\overline{\ze}} = \frac{h'}{D^2} \le( -1 +\tfrac{1}{2}(1-|\ze|^2)^2 (P_f)_{\overline{\ze}} \ri)
$$
and 
$$
\Phi_\ze = \frac{(1-|\ze|^2) h''}{D} + \frac{h'}{D^2}\le( -\overline{\ze}^2 +\tfrac{1}{2}(1-|\ze|^2)^2 (P_f)_{\ze} \ri). 
$$
Making use of \eqref{F-zed} we find that 
\begin{align*}
-z^2 F_z = \, &  \overline{\om h'} +\frac{h'}{D^2} \le( -1 +\tfrac{1}{2}(1-|\ze|^2)^2 (P_f)_{\overline{\ze}} \ri) + (1-|\ze|^2) \overline{\le(\frac{\om' h'}{D}\ri)} \\
 & + (1-|\ze|^2) \overline{\le(\frac{\om h''}{D}\ri)} + \overline{\le(\frac{\om h'}{D^2}\ri)} \le( -\ze^2 +\tfrac{1}{2}(1-|\ze|^2)^2 \overline{(P_f)_{\ze}} \ri). 
\end{align*}
For $\la\in\partial\SD$ as in \eqref{lambda-T} we write 
\begin{align*}
-z^2 F_z = \,  \frac{h'}{D^2}& \Bigg[  \la \, \overline{\om D^2} -1 +\tfrac{1}{2}(1-|\ze|^2)^2 (P_f)_{\overline{\ze}} + (1-|\ze|^2) \la \, \overline{\om' D } \\
 &  + (1-|\ze|^2) \la \, \overline{ \om D Ph} + \la\, \overline{ \om } \le( -\ze^2 +\tfrac{1}{2}(1-|\ze|^2)^2 \overline{(P_f)_{\ze}} \ri) \Bigg]. 
\end{align*}
We compute 
$$
(P_f)_{\overline{\ze}} = - \frac{|\om'|^2}{(1-|\om|^2)^2}\, , 
$$
so that, in view of \eqref{def-hyp-derivative}, we have that $(1-|\ze|^2)^2 (P_f)_{\overline{\ze}}=-|\om^*|^2$. We substitute this into the above, as well as $(P_f)_{\ze}$ from \eqref{def-Sch-HM} and $Ph$ from \eqref{pre-Sch-HM}, in order to obtain
\begin{align*}
-z^2 F_z = \, &  \frac{h'}{D^2} \Bigg[  \la \, \overline{\om D^2} -1 -\tfrac{1}{2}|\om^*|^2 + (1-|\ze|^2) \la \, \overline{\om' D } \\
  + & (1-|\ze|^2) \la \, \overline{ \om D}  \le( \overline{P_f} + \frac{ \overline{\om'}\om }{1-|\om|^2}  \ri)  - \la\overline{ \om }\ze^2 + \tfrac{\la}{2}\overline{\om}(1-|\ze|^2)^2 \le( \overline{S_f} + \tfrac{1}{2} \overline{P_f}^2 \ri) \Bigg]. 
\end{align*}
Substituting $\om'$ from \eqref{def-hyp-derivative} and the expression for $D$ we get
$$
-z^2 F_z =  \frac{h'}{D^2} \bigg[  -1 +  \la \ze \, \overline{\om^*}  -\tfrac{1}{2}|\om^*|^2  -  \tfrac{\la}{2}(1-|\ze|^2) \overline{\om^* P_f} + \tfrac{\la}{2}(1-|\ze|^2)^2 \overline{\om S_f}  \bigg]. 
$$
On the other hand, in view of \eqref{F-zed-bar} we have that 
$$
-\overline{z}^2 F_{\overline{z}} = \frac{h'}{D^2}  \bigg[ D^2 + (1-|\ze|^2) D Ph - \overline{\ze}^2 +\tfrac{1}{2}(1-|\ze|^2)^2 (P_f)_{\ze}   - \la  \overline{\om}  +\tfrac{\la}{2} \overline{\om} (1-|\ze|^2)^2 \overline{(P_f)_{\overline{\ze}}} \bigg], 
$$
for which we work as before in order to obtain
$$
-\overline{z}^2 F_{\overline{z}} = \frac{h'}{D^2} \bigg[ \overline{\ze\om} \om^* - \la\overline{\om}  - \tfrac{\la}{2} \overline{\om} |\om^*|^2 - \tfrac{1}{2} (1-|\ze|^2) \overline{\om} \om^* P_f + \tfrac{1}{2} (1-|\ze|^2)^2 S_f \bigg]. 
$$
Hence, factoring out $\om^*$, for the Beltrami coefficient we have that 
\begin{align*}
\le(\frac{\overline{z}}{z}\ri)^2 \mu_F = & \, \frac{-\overline{z}^2 F_{\overline{z}}}{-z^2 F_z} \\
= & \, \frac{ - \la\overline{\om} + \tfrac{1}{2} (1-|\ze|^2)^2 S_f  +  \om^* \le( \overline{\ze\om}  - \tfrac{\la}{2} \overline{\om\om^*}   - \tfrac{1}{2} (1-|\ze|^2) \overline{\om}  P_f \ri) }{  -1 + \tfrac{\la}{2}(1-|\ze|^2)^2 \overline{\om S_f}  + \overline{\om^*} \le(  \la \ze  -\tfrac{1}{2}\om^* -  \tfrac{\la}{2}(1-|\ze|^2) \overline{ P_f} \ri)  }. 
\end{align*}
We set 
$$
\al = \la \overline{\om} \qquad \text{and} \qquad \beta = \tfrac{1}{2}(1-|\ze|^2)^2 S_f
$$
and note that $|\al| \leq d <1$ and $|\beta| \leq t/2$. Dividing both numerator and denominator by $1-\overline{\beta}\al$ we obtain 
$$
\le(\frac{\overline{z}}{z}\ri)^2 \mu_F = \frac{ \vp_\beta(\al) -  \om^* M }{ 1  - \overline{\om^*}  N  }, 
$$ 
where $\vp_\beta(\al)$ is as in \eqref{disk-auto}, 
$$
M = \frac{ \overline{\ze\om}  - \tfrac{\la}{2} \overline{\om\om^*}   - \tfrac{1}{2} (1-|\ze|^2) \overline{\om}  P_f }{1-\overline{\beta}\al }
$$
and 
$$
N = \frac{  \la \ze  -\tfrac{1}{2}\om^* -  \tfrac{\la}{2}(1-|\ze|^2) \overline{ P_f} }{1-\overline{\beta}\al }. 
$$
Now, clearly 
$$
 |1- \overline{\beta}\al|  \geq 1- |\al\beta| \geq 1- dt/2  
$$
and, also, in view of \eqref{pre-Sch-HM}, \eqref{Pom-lem} and \eqref{Sh-small} we have that 
$$
(1-|\ze|^2) | P_f | \leq  (1-|\ze|^2) |Ph| + |\om\om^*| \leq 2 +2 \sqrt{ 1+ \tfrac{t+C R_t}{2} }+ d R_t \leq k 
$$
for some constant $k>0$. Therefore, 
$$
|M| \leq \frac{d +  d R_t/2 +d k/2}{1- dt/2  } \leq K
$$
and
$$
|N| \leq \frac{1 + R_t/2 +k/2}{1- dt/2  } \leq L
$$
for constants $K,L>0$. By the Schwarz-Pick lemma we have that 
$$
|\vp_\beta(\al)|  \leq \frac{|\al| + |\beta|}{1+ |\al\beta|} \leq \frac{d + t/2 }{1+ dt /2 }  
$$
and, finally, for the Beltrami coefficient we have the bound 
$$
|\mu_F| \leq \frac{ \frac{d + t/2 }{1+ dt /2 } +KR_t }{ 1- LR_t }. 
$$
Seeing that this bound tends to $d$ as $t\to0$ the proof is complete.

\vskip.3cm
\emph{Acknowledgements}. 
The second author was supported by Fondecyt grant \linebreak \#1190756, Chile. The third author was supported by the Spanish research projects PID2019-106093GB-I00 and PID2019-106870GB-I00 (MICINN); and  \linebreak  Spanish Thematic Research Network MTM2017-90584-REDT.

\end{document}